\documentclass[12pt]{amsart}
\begin{document}
\numberwithin{equation}{section}

\def\1#1{\overline{#1}}
\def\2#1{\widetilde{#1}}
\def\3#1{\widehat{#1}}
\def\4#1{\mathbb{#1}}
\def\5#1{\frak{#1}}
\def\6#1{{\mathcal{#1}}}

\def\C{{\4C}}
\def\R{{\4R}}
\def\N{{\4N}}
\def\Z{{\4Z}}

\title{On different notions of homogeneity for CR-manifolds}
\author[D. Zaitsev]{Dmitri Zaitsev}
\dedicatory {Dedicated to M. Salah Baouendi on the occasion of his birthday}
\address{D. Zaitsev: School of Mathematics, Trinity College Dublin, Dublin 2, Ireland}
\email{zaitsev@maths.tcd.ie}
\thanks{The author was supported in part by the RCBS grant of the Trinity College Dublin
and by the Science Foundation Ireland grant 06/RFP/MAT018}

\begin{abstract}
We show that various notions of local homogeneity for CR-manifolds are equivalent.
In particular, if germs at any two points of a CR-manifold are CR-equivalent,
there exists a transitive local Lie group action by CR-automorphisms near every point.
\end{abstract}

\maketitle

\def\Label#1{\label{#1}}


\def\cn{{\C^n}}
\def\cnn{{\C^{n'}}}
\def\ocn{\2{\C^n}}
\def\ocnn{\2{\C^{n'}}}


\def\dist{{\rm dist}}
\def\const{{\rm const}}
\def\rk{{\rm rank\,}}
\def\id{{\sf id}}
\def\aut{{\sf aut}}
\def\Aut{{\sf Aut}}
\def\CR{{\rm CR}}
\def\GL{{\sf GL}}
\def\Re{{\sf Re}\,}
\def\Im{{\sf Im}\,}
\def\span{\text{\rm span}}

\def\codim{{\rm codim}}
\def\crd{\dim_{{\rm CR}}}
\def\crc{{\rm codim_{CR}}}

\def\phi{\varphi}
\def\eps{\varepsilon}
\def\d{\partial}
\def\a{\alpha}
\def\b{\beta}
\def\g{\gamma}
\def\G{\Gamma}
\def\D{\Delta}
\def\Om{\Omega}
\def\k{\kappa}
\def\l{\lambda}
\def\L{\Lambda}
\def\z{{\bar z}}
\def\w{{\bar w}}
\def\Z{{\1Z}}
\def\t{\tau}
\def\th{\theta}

\emergencystretch15pt
\frenchspacing

\newtheorem{Thm}{Theorem}[section]
\newtheorem{Cor}[Thm]{Corollary}
\newtheorem{Pro}[Thm]{Proposition}
\newtheorem{Lem}[Thm]{Lemma}

\theoremstyle{definition}\newtheorem{Def}[Thm]{Definition}

\theoremstyle{remark}
\newtheorem{Rem}[Thm]{Remark}
\newtheorem{Exa}[Thm]{Example}
\newtheorem{Exs}[Thm]{Examples}

\def\bl{\begin{Lem}}
\def\el{\end{Lem}}
\def\bp{\begin{Pro}}
\def\ep{\end{Pro}}
\def\bt{\begin{Thm}}
\def\et{\end{Thm}}
\def\bc{\begin{Cor}}
\def\ec{\end{Cor}}
\def\bd{\begin{Def}}
\def\ed{\end{Def}}
\def\br{\begin{Rem}}
\def\er{\end{Rem}}
\def\be{\begin{Exa}}
\def\ee{\end{Exa}}
\def\bpf{\begin{proof}}
\def\epf{\end{proof}}
\def\ben{\begin{enumerate}}
\def\een{\end{enumerate}}
\def\beq{\begin{equation}}
\def\eeq{\end{equation}}

\section{Introduction}

The purpose of this paper is to show that various notions of local homogeneity
for real-analytic CR-manifolds are in fact equivalent.
The case of real-analytic hypersurfaces $M$ in $\C^2$ has
been considered by A.~V.~Loboda in \cite{L98}, where the equivalence of two different notions is shown,
namely biholomorphic equivalence of germs of $M$ at any two points
and the existence of a transitive local Lie group action via biholomorphisms near every point of $M$.
The proof is based on a refined Chern-Moser normal form \cite{CM}
and convergence radius estimates due to V.~K.~Beloshapka and A.~G.~Vitushkin \cite{BV}.
In this paper we extend this result to arbitrary real-analytic CR-manifolds,
for which no such normal form is available in general.
We also propose weaker homogeneity conditions based on the notion of ``$k$-equivalence''
introduced in \cite{BRZ} and show that they still lead to an equivalent notion of local homogeneity.

These results appear to be in sharp contrast with the fact that different
nonequivalent notions exist for {\em global} homogeneity.
In fact, W.~Kaup \cite{K67} constructed an example of a domain $D\subset\C^2$,
which is homogeneous in the sense that any two points
are mapped into each other by a (global) biholomorphic automorphism of $D$
but no (finite-dimensional) Lie group acts transitively on $D$
via biholomorphic automorphisms.

We now briefly recall the necessary definitions to state our results.
The reader is referred e.g.\ to the book \cite{BERbook} for further details and related facts.
An {\em (abstract) CR-manifold} is a real manifold $M$ together with a 
formally integrable distribution $\6V$ of its complexified tangent space $\C TM$
satisfying $\6V\cap \1{\6V}=0$, called the {\em CR-structure}
(here $\1{\6V}$ denotes the complex conjugate subbundle).
A {\em CR-map} between two CR-manifolds $M$ and $M'$ with CR-structures $\6V$
and $\6V'$ is any map $h\colon M\to M'$ with $h_*\6V\subset V'$,
a CR-diffeomorphism is any diffeomorphism, which  is CR together with its inverse,
and a CR-automorphism is CR-diffeomorphisms from a manifold into itself.
{\em All CR-manifolds and CR-maps in this paper will be assumed to be real-analytic.}
This is motivated by our primary interest in homogeneous CR-manifolds and
the fact that any CR-manifold admitting a transitive Lie group
action by CR-automorphisms (or even a transitive local Lie group action, see below)
is automatically real-analytic. 

It is well-known that a real-analytic
CR-manifold is locally embeddable into $\C^N$ with suitable $N$
such that its CR-structure is induced by the complex structure of $\C^N$
(which is a special case of a CR-structure with $\6V= T^{(0,1)}\C^N$).
This allows to pass from intrinsic to extrinsic point of view and vice versa,
which we shall frequently do here.

Two germs $(M,p)$ and $(M',p')$ of CR-manifolds are said to be {\em CR-equivalent}
if there exists a CR-diffeomorphism between open neighborhoods of $p$ and $M$
and of $p'$ in $M'$ sending $p$ into $p'$. A weaker notion is that of a 
{\em formal CR-equivalence}, where $(M,p)$ and $(M',p')$ are said to be 
{\em formally CR-equivalent} if there exists an invertible formal power series map
(in some and hence any local real-analytic coordinates on $M$ and $M'$)
which sends $\6V$ into $\6V'$ in the formal sense.
Yet more generally, $(M,p)$ and $(M',p')$ are said to be {\em $k$-equivalent},
where $k>1$ is any integer, if there exists an invertible real-analytic map $h$
between open neighborhoods of $p$ and $M$
and of $p'$ in $M'$ sending $p$ into $p'$
and sending $\6V$ into $\6V'$ ``up to order $k$''.
The latter means that given a local frame $e_1(x),\ldots,e_n(x)\in \6V_x$, $x\in M$,
one can find a corresponding frame $e'_1(x'),\ldots,e'_n(x')\in \6V_{x'}$, ${x'}\in M$,
of $\6V'$ such that $h_*e_j(x) = e'_j(h(x)) + O(|x|^k)$,
where $x\in \R^{\dim M}$ is any local coordinate system vanishing at $p$.
By a result of M. S. Baouendi, L. P. Rothschild and the author \cite[Corollary~1.2]{BRZ},
the notions of being CR-equivalent, formally CR-equivalent and $k$-equivalent for all $k$
are equivalent for germs of CR-manifolds at their points in {\em general position}
(see Theorem~\ref{decomp}~(iv)).
On the other hand, a similar fact does not hold for more general
real-analytic submanifolds in $\C^N$ in view of an example by J. Moser and S. Webster \cite{MW}.
It is an open question whether the same conclusion holds for arbitrary CR-manifolds.

Another type of notion of local homogeneity is based on local Lie group actions.
As customary we always assume a Lie group to have at most countably many connected components.
Recall that a (real-analytic) {\em local action} of a Lie group $G$ with unit $e$ on a manifold $M$
is a neighborhood $O$ of $\{e\}\times M$ and a real-analytic map 
$\phi\colon O\to M$, $(g,x)\mapsto g\cdot x$,
satisfying $e\cdot x=x$ and $(g_1 g_2) \cdot x = g_1 \cdot (g_2 \cdot x)$
whenever both sides are defined (see \cite{P} for further details on local group actions).
A local Lie group action is said to be {\em transitive} if for every $p,q\in M$, 
there exists a finite sequence $g_1,\ldots,g_s\in G$ such 
that all expressions $R_j:=g_j\cdot (g_{j-1}\cdot\ldots (g_1\cdot p) \ldots)$ are defined
for $1\le j\le s$ and $R_s=q$.
It is easy to see that if $M$ is connected, 
a local Lie group action $\phi$ is transitive
if and only if the differential $\phi_*$ sends $T_eG\subset T_{(e,p)} (G\times M)$ onto $T_pM$
for every $p\in M$.

Another, a priori weaker notion is based on the following generalization
of a transitive local Lie group action that we state in a local form for germs:
\bd\Label{tr-family}
We say that a germ of CR-manifold $(M,p)$ admits a 
{\em transitive family of local CR-automorphisms}
if there exists a germ of a CR-map $\phi\colon (M,p)\times (\R^{\dim M},0) \to (M,p)$
with $\phi_*(T_pM) = \phi_*(T_0\R^{\dim M}) = T_pM$, where $\phi_*$ is taken at $(p,0)$.
\ed

Finally we make use of (real-analytic) {\em infinitesimal CR-automorphisms} of $M$,
which are (real-analytic) real vector fields on $M$ whose local flows
are $1$-parameter families of CR-automorphisms
(see e.g.\ \cite[\S12.4]{BERbook} for more details).

We can now state our global result:

\bt\Label{global}
Let $M$ be a connected (real-analytic) CR-manifold.
Then the following are equivalent:
\ben
\item[(i)] for every $p,q\in M$, the germs $(M,p)$ and $(M,q)$ are $k$-equivalent for every $k$;
\item[(ii)] for every $p,q\in M$, the germs $(M,p)$ and $(M,q)$ are formally CR-equivalent;
\item[(iii)] for every $p,q\in M$, the germs $(M,p)$ and $(M,q)$ are CR-equivalent;
\item[(iv)] for every $p\in M$, the germ $(M,p)$ admits a transitive family of local CR-automorphisms;
\item[(v)] for every $p\in M$, the germs of all infinitesimal CR-automorphisms of $(M,p)$ span the tangent space $T_pM$;
\item[(vi)] for every $p\in M$, there exists a finite-dimensional Lie algebra of germs of infinitesimal CR-automorphisms of $(M,p)$ that spans the tangent space $T_pM$;
\item[(vii)] there exists a Lie group $G$ and, for every $p\in M$, a transitive local action of $G$
 by CR-automorphisms on an open neighborhood of $p$ in $M$.
\een
\et

In our second main result we refine Theorem~\ref{global}
stating all local homogeneity conditions for a germ of a CR-manifold $(M,p)$,
where the homogeneity means that some representative of the germ is locally homogeneous.
It turns out that the weakest condition (i) in Theorem~\ref{global}
can be here further weakened by requiring that only the germs of $M$
at sufficiently many points are equivalent rather than all germs.
A more precise definition is as follows.

\bd\Label{weak-orbit}
Let $M$ be a real-analytic CR-manifold and $p\in M$ be an arbitrary point.
The {\em weak equivalence orbit} of $p$ in $M$ is the set of all $q\in M$
such that the germs $(M,q)$ and $(M,p)$ are $k$-equivalent for all $k$.
We say that the germ $(M,p)$ {\em satisfies condition $(*)$} if 
for any open neighborhood $U$ of $p$ in $M$,
the weak equivalence orbit of $p$ in $U$ is not contained in
a real-analytic submanifold of $U$ of smaller dimension.
\ed

We can now state our local result.

\bt\Label{local}
Let $(M,p)$ be a germ of a real-analytic CR-manifold,
where we write $M$ for any representative.
Then the following are equivalent:
\ben
\item[(i)] $(M,p)$ satisfies condition $(*)$;
\item[(ii)] the weak equivalence orbit of $p$ in $M$ contains an open neighborhood of $p$ in $M$;
\item[(iii)] for every $q\in M$ sufficiently close to $p$, 
the germs $(M,q)$ and $(M,p)$ are formally CR-equivalent;
\item[(iv)] for every $q\in M$ sufficiently close to $p$, 
the germs $(M,q)$ and $(M,p)$ are CR-equivalent;
\item[(v)] $(M,p)$ admits a transitive family of local CR-automorphisms;
\item[(vi)] the germs of all infinitesimal CR-automorphisms of $(M,p)$ span the tangent space $T_pM$;
\item[(vii)] there exists a finite-dimensional Lie algebra of germs of infinitesimal CR-automorphisms of $(M,p)$ that spans the tangent space $T_pM$;
\item[(viii)] there exists a Lie group $G$ and a transitive local action of $G$
 by CR-automorphisms on an open neighborhood of $p$ in $M$.
\een
\et

The proofs of Theorems~\ref{global} and \ref{local} are given in \S\ref{proofs}.
In \S\ref{struc} we state basic structure results for general CR-manifolds and their maps
that play crucial role in the proofs.
In \S\ref{definable} we recall few definitions and facts about 
sets definable in terms of certain rings of functions,
also needed for the proofs.
In \S\ref{technical} we prove a proposition that represents
the main technical core of the proofs of Theorems~\ref{global} and \ref{local}.

We conclude by mentioning that (locally) homogeneous CR-manifold
can be described in purely algebraic terms, e.g.\ in terms of the
so-called ``CR-algebras'' considered in \cite{MN}, see also \cite{F}
for a local description.

\section{Structure results for CR-manifolds and jet parametrization of CR-diffeomorphisms}\Label{struc}

We recall here some basic definition and structure results from \cite{BRZ} 
for real-analytic CR-manifolds.
We first note that any real-analytic CR-manifold can be locally
embedded as a real-analytic generic submanifold into $\C^N$ for suitable $N$
(see e.g.\ \cite[Chapter II]{BERbook}).
(Recall that a real submanifold $M\subset\C^N$ is {\em generic} if $T_pM + iT_pM = T_p\C^N$
for every $p\in M$.)
We thus give the extrinsic definitions for embedded generic submanifolds of $\C^N$
following \cite{BRZ} that will suffice for our purposes
(see e.g. \cite[Chapter XI]{BERbook} for an intrinsic approach).
Let $\rho(Z,\1Z)=(\rho^1(Z,\1Z),\ldots,\rho^d(Z,\1Z))$ 
be a vector-valued local defining function of $M$ near a point $p$,
i.e.\ with the rank of $\frac{\d\rho}{\d Z}$ being equal to the codimension of $M$.
Recall that a $(0,1)$ vector field on $M$ is any vector field of the form 
$L=\sum a_j(Z,\1Z)\frac{\d}{\d\1Z_j}$ with $(L\rho)(Z,\1Z)\equiv 0$ on $M$.
In our case when $M$ is real-analytic,
it will be sufficient to consider only real-analytic vector fields.

Following \cite[\S2.3]{BRZ}, consider the vector subspace
\begin{equation}\Label{rp}
E(p):=\span_\C \left\{(L_1\ldots L_s\rho^j_Z)(p,\1p) :  1\le j\le d;
0\le s <\infty; \right\} \subset \C^N,
\end{equation}
where $L_1,\ldots,L_s$ run through all collections of $(0,1)$ vector fields
and $\rho_Z^j(Z,\1Z)\in\C^N$ denotes the
complex gradient of $\rho^j$ with respect to $Z$.
The number $r_2(p):=N-\dim_\C E(p)$ 
is said to be the {\em degeneracy} of $M$ at $p$
and $M$ is said to be of {\em minimum degeneracy}
at a point $p_0\in M$ if $p_0$ is a local minimum of the integer function $p\mapsto r_2(p)$.
Recall that $(M,p)$ is {\em finitely nondegenerate} if and only if $r_2(p)=0$
(i.e.\ $E(p)=\C^N$)
and is {\em $l$-nondegenerate} if $l$ is the smallest integer such 
that $\C^N$ is spanned by the vectors $L_1\ldots L_s\rho^j_Z)(p,\1p)$ with $s\le l$.

Similarly consider the vector subspace $\frak g_M(p)$ of the complexified tangent space
$\C T_pM$ generated by the values at $p$
of all $(0,1)$ vector fields, their conjugates and all finite order commutators
involving $(0,1)$ vector fields and their conjugates.
The corresponding number $r_3(p):=\dim_\R M -dim_\C \frak g_M(p)$
is said to be the {\em orbit codimension} of $M$ at $p$
and $M$ is said to be of {\em minimum orbit codimension}
at a point $p_0\in M$ if $p_0$ is a local minimum of the function $p\mapsto r_3(p)$.
Recall that $(M,p)$ is of {\em finite type}
(in the sense of Kohn and  and Bloom-Graham) if and only if $r_3(p)=0$
(i.e.\ $\frak g_M(p)=\C T_pM$).

The following theorem summarizes some of the results
by M. S. Baouendi, L. P. Rothschild \cite{BRZ}
that will be crucial for the proofs of both Theorems~\ref{global} and \ref{local}.

\bt\Label{decomp}
Let $M$ be a connected real-analytic CR-manifold
and $V\subset M$ be the subset of all points $p\in M$
such that $M$ is either not of minimum degeneracy or not of minimum orbit codimension at $p$.
Then $V$ is a (closed) proper real-analytic subvariety of $M$ and there exist
nonnegative integers $N_1,N_2,N_3$ and, for every $p\in M\setminus V$, 
a generic real-analytic submanifold $M'\subset \C^{N_1}\times \R^{N_2}\subset \C^{N_1+N_2}$ 
passing through $0$ such that the following hold:
\ben
\item[(i)] $(M,p)$ is CR-equivalent to $(M'\times \C^{N_3},0)$;
\item[(ii)] $(M',0)$ is finitely nondegenerate;
\item[(iii)] for every $u\in \R^{N_2}$ near $0$, one has
$(0,u)\in M'$ and $M'\cap (\C^{N_1}\times\{u\})$ 
is a CR-manifold of finite type at $(0,u)$;
\item[(iv)] if $(\2M,q)$ is another germ of a real-analytic CR-manifold
such that $(M,p)$ and $(\2M,q)$ are $k$-equivalent for any $k$, then they are also CR-equivalent.
\een
\et

In fact, for $u\in \R^{N_2}$ near $0$, 
the slice $(M'\cap (\C^{N_1}\times\{u\})$ represents the so-called 
{\em local CR-orbit} of $M'$ at $(0,u)$.
Recall that the local {\em CR orbit} of
a point $q\in M'$ is the germ at $q$ of a
(real-analytic) submanifold of $M'$ through $q$ of smallest
possible dimension to which all the
$(0,1)$ vector fields on $M$ are tangent.
(The existence and uniqueness of a local CR-orbit
is a consequence of a theorem of Nagano~\cite{N66},
see also \cite[\S3.1]{BERbook}.)
Note that in general $M'$ cannot be locally written as a product 
of its CR-orbit and $\R^{N_2}$ since different CR-orbits may not
be CR-equivalent (see \cite[\S2]{BRZproc} for an example).

\br
The integers $N_1,N_2,N_3$ are uniquely determined by $M$,
where $N_2$ is the minimum degeneracy
and $N_3$ the minimum orbit codimension of $M$,
see \cite{BRZ}.
\er

We shall also need the following result from \cite{BRZ}
(see also \cite{BRZproc}) describing the behavior
of CR-equivalences with respect to the decomposition
provided by Theorem~\ref{decomp}.

\bt\Label{decomp-map}
Let $M,M'\subset \C^{N_1}\times \R^{N_2}$ be generic real-analytic
submanifolds of the same dimension passing through $0$, 
both satisfying {\rm(ii)} and {\rm(iii)} of Theorem~\ref{decomp},
i.e.\  such that both $(M,0)$ and $(M',0)$ are finitely nondegenerate
and for every $u\in \R^{N_2}$ near $0$, one has
$(0,u)\in M\cap M'$ and  $M\cap (\C^{N_1}\times\{u\})$ and $M'\cap (\C^{N_1}\times\{u\})$ 
are both of finite type at $(0,u)$.
Let $H\colon (\C^{N_1}\times\C^{N_2}\times\C^{N_3},0)\to (\C^{N_1}\times\C^{N_2}\times\C^{N_3},0)$ 
be a germ of a biholomorphic map fixing $0$ and sending 
$M\times\C^{N_3}\times\{0\}$ into $M'\times\C^{N_3}\times\{0\}$.
Then $H$ is of the form
\begin{equation}\Label{form}
H(Z_1,Z_2,Z_3) = (H_1(Z_1,Z_2),H_2(Z_2),H_3(Z_1,Z_2,Z_3)),
\quad (Z_1,Z_2,Z_3)\in \C^{N_1}\times\C^{N_2}\times\C^{N_3},
\end{equation}
where $H_2$ is a local biholomorphic map of $\C^{N_2}$ preserving $\R^{N_2}$
and for $u\in\R^{N_2}$ near $0$, 
$H_1(\cdot,u)$ a local biholomorphic map of $\C^{N_1}$
sending $(M\cap (\C^{N_1}\times\{u\}), (0,u))$ into 
$(M'\cap (\C^{N_1}\times\{H_2(u)\}), (0,H_2(u))$
(both regarded as generic submanifolds of $\C^{N_1}$).
\et

Our next main ingredient is a parametrization
result from \cite{BERrat} for local biholomorphisms between generic manifolds with parameters.
(Further parametrization results of this kind can be found
in \cite{Ejets,ELZ,KZ,LM,LMZ}.)
Here we consider a {\em real-analytic family of generic submanifolds of $\C^N$},
which is a collection $\{M_x\}$ of generic submanifolds of $\C^N$ 
with parameter $x$ from another real-analytic manifold $X$ such that for 
every $x_0\in X$ and $p\in M_{x_0}$,
all manifolds $M_x$ near $p$ with $x\in X$ near $x_0$
can be defined by a family of defining functions $\rho(Z,\bar Z,x)$,
which is real-analytic in all its arguments.
We also write $J^k_{0,0}(\C^N,\C^N)$
for the space of all $k$-jets of holomorphic maps from $\C^N$ into itself with
both source and target being $0$.

\bt\Label{parametrization}
Let $M_x$, $x\in X$, and $M'_{x'}$, $x'\in X'$, 
be real-analytic families of generic submanifolds through $0$ in $\C^N$ of codimension $d$.
Assume that, for some fixed points $x_0\in X$ and $x'_0\in X_0$,
\ben
\item[(i)] $M_{x_0}$ is of finite type at $0$;
\item[(ii)] $M'_{x_0}$ is $l$-nondegenerate at $0$ for some $l$.
\een
Set $k:=l(d+1)$.
Then for every invertible jet $\L_0\in J^k_{0,0}(\C^N,\C^N)$, there exist open neighborhoods
$\Omega'$ of $0$ in $\C^N$, $\Omega''$ of $\L_0$  in  $J^k_{0,0}(\C^N,\C^N)$,  
$U$ of $x_0$ in $X$ and $U'$ of $x'_0$ in $X'$,
and a real-analytic map 
$\Psi\colon\Omega'\times \Omega'' \times U\times U'\to \C^N$
such that the identity
\beq\Label{jps}
H(Z) = \Psi(Z,j^{k}_0 H, x, x')
\eeq
holds for any $x\in U$, $x'\in U'$, 
any local biholomorphism $H$ of $\C^N$ fixing $0$ and sending $M_x$ into $M'_{x'}$
and any $Z\in \Omega'$ sufficiently close to $0$.
\et

Finally, in the setting of Theorem~\ref{parametrization}, 
it will be important to describe the sets of those jets
that actually arise as jets of local biholomorphisms between $(M_x,0)$ and $(M'_{x'},0)$.
We shall make use of the following result,
similar to \cite[Theorem~5.2.9]{BERrat} whose
proof can be obtained by repeating the corresponding arguments in \cite{BERrat}:

\bt\Label{anal-jets}
Under the assumptions of Theorem~\ref{parametrization},
there exist open neighborhoods $U$ of $x_0$ in $X$ and $U'$ of $x'_0$ in $X'$,
and finite sets of polynomials $a_j(\L,\1\L,x,x')$ and $b_s(\L,x,x')$ in 
$(\L,\1\L)\in J^k_{0,0}(\C^N,\C^N)\times \1{J^k_{0,0}(\C^N,\C^N)}$
with real-analytic coefficients in $(x,x')\in U\times U'$
such that the set of all $(\L,x,x')\in J^k_{0,0}(\C^N,\C^N)\times U\times U'$,
for which there exists a local biholomorphism $H$ of $\C^N$ fixing $0$
and sending $M_x$ into $M'_{x'}$ with $j^k_0 H = \L$, is given by
\begin{equation}\Label{set}
\{ a_j(\L,\1\L,x,x') = 0 \text{ for all } j\} \setminus 
\{ b_s(\L,x,x') = 0 \text{ for all } k\}.
\end{equation}
\et

\section{Definable and semianalytic sets}\Label{definable}

Here we collect some basic definitions and properties of
sets definable over rings, in particular, of semianalytic sets.
The readers is referred e.g.\ to \cite{BM}
and the extensive literature cited there for proofs and further related facts.

Let $\6R$ be a ring of real-valued functions on a set $E$.
A subset $A\subset E$ is said to be {\em definable} over $\6R$
if $A$ can be written as $\cup_{j=1}^s \cap_{k=1}^r A_{jk}$,
where each $A_{jk}$ is either $\{f_{jk}=0\}$ or $\{f_{jk}>0\}$
with $f_{jk}\in \6R$.
In particular, a subset $A$ in a real-analytic manifold $M$
is called {\em semianalytic} if every point $p\in M$
has an open neighborhood $U$ such that $A\cap U$ is definable over 
the ring of all real-analytic functions on $U$.
It is elementary to see that any real-analytic subset is always semianalytic
and that finite unions, intersections and complements
of semianalytic sets are again semianalytic.

The following is a fundamental structure theorem for semianalytic sets:

\bt\Label{semianal-struct}
Every semianalytic set $A\subset M$ admits a stratification
into a locally finite disjoint union of real-analytic submanifolds $A_j$ of $M$,
each being a semianalytic subset of $M$, and satisfying the ``frontier condition'':
if $A_j\cap \1{A_k}\ne \emptyset$, then $A_j\subset \1{A_k}$ and $\dim A_j < \dim A_k$.
\et

As a consequence, the Hausdorff dimension $\dim A$ equals to the maximum stratum dimension.
We shall use the following Lojaciewicz's version of the Tarski-Seidenberg theorem
(see e.g.\ \cite[\S2]{BM}):

\bt\Label{proj}
Let $\6R$ be a ring of functions on a set $E$
and $\6R[x_1,\ldots,x_k]$ be the corresponding polynomial ring on $E\times \R^k$.
Denote by $\pi\colon E\times \R^k\to E$ the canonical projection.
Then, if $A\subset E\times \R^k$ is definable over $\6R[x_1,\ldots,x_k]$,
its projection $\pi(A)\subset E$ is definable over $\6R$.

In particular, if $M$ is a real-analytic manifold and $A\subset M\times \R^k$
is definable over the ring of polynomials in $(x_1,\ldots,x_k)\in\R^k$ with 
real-analytic coefficients in $M$, then its projection $\pi(A)\subset M$
is semianalytic.
\et

Note that it is essential in Theorem~\ref{proj} that $A$
is definable over the ring of polynomials with real-analytic coefficients in $M$
rather than the ring of all real-analytic functions on $M\times \R^k$,
for which the corresponding conclusion would fail
(see e.g.\ \cite[\S2]{BM} for an example).

\section{Weak equivalence orbits and their properties}\Label{technical}

Here we consider weak equivalence orbits and condition $(*)$ as
defined in Definition~\ref{weak-orbit} and 
obtain its implications that will be crucial
for the proofs of Theorems~\ref{global} and \ref{local}.
As before $M$ denotes a connected real-analytic CR-manifold and $p\in M$ its arbitrary point.

\bl\Label{not-contained}
Let $(M,p)$ satisfy condition $(*)$.
Then the weak equivalence orbit of $p$ in $M$
is not contained in any semianalytic subset $A\subset M$ 
with $\dim A<\dim M$.
\el

\bpf
Without loss of generality, $M$ is connected.
Assume, by contradiction, that the weak equivalence orbit $O$ of $p$ in $M$
is contained in a semianalytic subset $A\subset M$ of lower dimension.
Fix a stratification of $A$
into a locally finite disjoint union of real-analytic submanifolds $A_j$, 
that exists due to Theorem~\ref{semianal-struct}.
Let $m$, $0\le m <\dim M$, be the minimum integer such that $O$ is contained
in the union $\2A$ of all strata of dimension not greater than $m$.
Then there exists a point $q\in O$ which is contained in a stratum $A_j$ of dimension precisely $m$.
Now the ``frontier condition'' in Theorem~\ref{semianal-struct} implies
that $q$ is not contained in the closure of any stratum $A_k$ with $\dim A_k\le \dim A_j = m$.
Hence, by the choice of $m$, there exists an open neighborhood $\Omega$ of $q$ in $M$
such that $O\cap \Omega\subset A_j$.
Finally, by the definition of the weak equivalence orbit, the germs $(M,q)$ and $(M,p)$
are $k$-equivalent for any $k$.
Hence they are also CR-equivalent in view of Theorem~\ref{decomp} (iv).
Let $\phi\colon U\to V$ be any CR-equivalence between open neighborhoods $U$
and $V$ of $q$ and $p$ respectively. Without loss of generality, $U\subset \Omega$.
Then $\phi$ sends $O\cap U$ onto $O\cap V$ and therefore $O\cap V$
is contained in the low dimensional submanifold $\phi(A_j\cap U)$ of $V$,
which is a contradiction with condition $(*)$. The proof is complete.
\epf

\bp\Label{local'}
Let $(M,p)$ satisfy condition $(*)$.
Then there exist integers $N_1,N_2,N_3$ and a generic real-analytic submanifold
$\2M\subset \C^{N_1}$ passing through $0$ such that the following hold:
\ben
\item[(i)] $(M,p)$ is CR-equivalent to $(\2M\times \C^{N_2}\times \R^{N_3},0)$;
\item[(ii)] $(\2M,0)$ is finitely nondegenerate and of finite type;
\item[(iii)] $(\2M,0)$ admits a transitive family of local CR-automorphisms.
\een
\ep

\bpf
Let $V\subset M$ be the proper real-analytic subvariety
considered in Theorem~\ref{decomp}.
Since $V$ is also a semianalytic subset of $M$ of a smaller dimension,
Lemma~\ref{not-contained} implies that the weak equivalence orbit $O$ of $p$ in $M$
is not contained in $V$. Let $q\in O\setminus V$ be any point.
Then $(M,q)$ is CR-equivalent to a germ $(M'\times \C^{N_3},0)$ as in Theorem~\ref{decomp}.
But since $q\in O$, the germs $(M,q)$ and $(M,p)$ are $k$-equivalent for any $k$
and therefore also CR-equivalent by Theorem~\ref{decomp} (iv).
Hence also $(M,p)$ is CR-equivalent to $(M'\times \C^{N_3},0)$.
Without loss of generality, we may assume $(M,p)=(M'\times \C^{N_3},0)$.
Since $(M,p)$ is assumed to satisfy condition $(*)$,
$(M'\times\C^{N_3},0)$ also does and
hence also $(M',0)$ satisfies condition $(*)$ in view of Theorem~\ref{decomp-map}.

We next consider for every $(q,u)\in M'\subset \C^{N_1}\times \R^{N_2}$,
the submanifold 
\begin{equation}\Label{constr}
\2M_{(q,u)}:= \{Z_1-q\in \C^{N_1} : (Z_1,u)\in M'\},
\end{equation}
passing through $0$.
It follows from our construction that $\2M_{(q,u)}$, $(q,u)\in M'$,
is a real-analytic family of generic submanifolds through $0$ in $\C^{N_1}$
and that $\2M_{(0,0)}$ is finitely nondegenerate and of finite type.
Hence we can apply Theorem~\ref{anal-jets}. 
As its consequence, we conclude that there exist an open neighborhood $U'$ of $(0,0)$ in $M'$
such that the set $A$ of all $(\L,x')\in J^k_{0,0}(\C^{N_1},\C^{N_1})\times U'$,
for which there exists a local biholomorphism of $\C^{N_1}$ 
sending $(\2M_{(0,0)},0)$ into $(\2M_{x'},0)$ with $j^k_0 H =\L$,
is definable (in the sense of \S\ref{definable}) over the ring of polynomials in $(\L,\1\L)$
with real-analytic coefficients in $x'\in U'$.
Here as in Theorem~\ref{anal-jets} we set $k:=l(d+1)$,
where $d$ is the codimension of $\2M_{(0,0)}$ in $\C^{N_1}$
and $l$ is such that $\2M_{(0,0)}$ is $l$-nondegenerate at $0$.

We now consider the natural projection $\pi\colon J^k_{0,0}(\C^{N_1},\C^{N_1})\times U' \to U'$.
Then $\pi(A)$ is a semianalytic subset of $U'$ by Theorem~\ref{proj}.
We claim that $\pi(A)$ contains the weak equivalence orbit of $0$ in $M'$.
Indeed, let $x'\in U'$ be in the orbit.
Note that, by our construction, $M'$ is both of minimum degeneracy
(in fact finitely nondegenerate) and of the minimum orbit codimension $d$ at $0$.
Then, in view of  Theorem~\ref{decomp} (iv), 
there exists a CR-equivalence $H'$ between germs $(M',0)$ and $(M',x')$,
which extends to a biholomorphic map of $\C^{N_1}\times\C^{N_2}$
sending $(M',0)$ into $(M',x')$ (see e.g.\ \cite[Corollary~1.7.13]{BERbook}
for the latter fact).
By Theorem~\ref{decomp-map}, $H$ is of the form \eqref{form} (with $N_3=0$).
Furthermore, it follows from the property of the component $H_1$ in
the decomposition \eqref{form} and our construction \eqref{constr} that 
$\2H(Z_1):= H(Z_1)-q$ is a local biholomorphism of $\C^{N_1}$ 
sending $(\2M_{(0,0)},0)$ into $(\2M_{x'},0)$, where $x'=(q,u)$.
But then $(j^k_0 \2H,x')\in A$ and hence $x'\in \pi(A)$, proving our claim.

We can now make use of Lemma~\ref{not-contained}
and conclude that the semianalytic subset $\pi(A)\subset U'$
must have the top dimension $\dim U'$.
Equivalently, $\pi(A)$ has a nonempty interior in $U'$.
Furthermore, the set $A\in J^k_{0,0}(\C^{N_1},\C^{N_1})\times U'$
is also semianalytic and hence admits itself a stratification 
in the sense of Theorem~\ref{semianal-struct}.
It follows that, in order for $\pi(A)$ to have a nonempty interior in $U'$,
there must exist a stratum $A_j$ of $A$ such that $\pi|_{A_j}\colon A_j\to U'$
is a submersion at some point of $A_j$.
By the implicit function theorem, there exists an open set $\Omega\subset U'$
and a real-analytic map $\nu\colon \Omega \to J^k_{0,0}(\C^{N_1},\C^{N_1})$
with $(\nu(x'),x')\in A_j\subset A$ for $x'\in \Omega$.

We next apply Theorem~\ref{parametrization} giving a parametrization $\Psi$
of local biholomorphisms $H$ sending $(\2M_{x},0)$ into $(\2M_{x'},0)$,
where we set $x_0:=0$, pick arbitrary $x'_0\in \Omega$ and set $\L_0:= \nu(x'_0)$.
Since for $x'\in\Omega$, we have $(\nu(x'),x')\in A$,
there exists a local biholomorphism $H_{x'}$ of $\C^{N_1}$
sending $(\2M_{(0,0)},0)$ into $(\2M_{x'},0)$,
which must therefore be given by the formula
$$H_{x'}(Z_1)=\Psi(Z_1,\nu(x'),0,x').$$
Then, in view of \eqref{constr}, for $x'=(q,u)\in M'\subset \C^{N_1}\times \R^{N_2}$
close to $x'_0$, the map 
\begin{equation}\Label{map-fam}
Z_1\mapsto\Psi(Z_1,\nu(x'),0,x') + q
\end{equation}
defines a local biholomorphism of $\C^{N_1}$ sending $(\2M_{(0,0)},0)$
into $(M'\cap (\C^{N_1}\times\{u\}),x')$, the latter being regarded as a submanifold of $\C^{N_1}$.

Let $S\subset M'$ be any real-analytic submanifold through $x'_0$ satisfying 
$$T_{x'_0} M' = T_{x'_0} \big(M'\cap(\C^{N_1}\times\{u\})\big) \oplus T_{x'_0} S.$$
Note that near $x'_0$, $S$ is automatically totally real and its projection to $\R^{N_2}$
defines a local diffeomorphism at $x'_0$ between $S$ and $\R^{N_2}$. Then the map 
$$(Z_1,x')\in \C^{N_1}\times S \mapsto \big(\Psi(Z_1,\nu(x'),0,x') + q,u\big)\in \C^{N_1}\times \R^{N_2},$$
where we keep the notation $x'=(q,u)$ as before,
defines a local CR-equivalence between $(\2M_{(0,0)},0)\times (S,x'_0)$ and $(M',x'_0)$. 
Since $S$ is totally real of dimension $N_2$, we conclude that
$(M',x'_0)$ is CR-equivalent to $(\2M_{(0,0)}\times \R^{N_2},(0,0))$.
Furthermore, by our construction, $(\nu(x'_0),x'_0)\in A$
implying that $(M',0)$ is CR-equivalent to $(M',x'_0)$ and therefore 
to $(\2M_{(0,0)},0)\times (\R^{N_2},0)$. This shows (i) and (ii)
with $\2M:= \2M_{(0,0)}$.

To show (iii), consider the family of local diffeomorphisms 
$\Phi_{x'}(Z_1):=\Psi(Z_1,\nu(x'),0,x') + q$.
We write $x'_0=(q_0,u_0)$ and let $x'=(q,u_0)\in M'\cap (\C^{N_1}\times \{u_0\})$.
Then, after a local identification of $(M'\cap (\C^{N_1}\times \{u_0\},x'_0))$
with $(\R^m,0)$, where $m$ is the corresponding dimension, the map 
$$(Z_1,q)\mapsto \Phi_{(q,u_0)}\circ \Phi_{(q_0,u_0)}^{-1}(Z_1)$$
defines a transitive family of local CR-automorphisms for $(M'\cap (\C^{N_1}\times \{u_0\},x'_0))$.
Finally, since $(M'\cap (\C^{N_1}\times \{u_0\},x'_0))$ is CR-equivalent to $(\2M,0)$,
the latter also admits a transitive family of local CR-automorphisms as desired.
The proof is complete.
\epf

\section{Proofs of Theorems~\ref{global} and \ref{local}}\Label{proofs}

We begin with Theorem~\ref{local}.
The implications (viii) $\Rightarrow$ (vii) $\Rightarrow$ (vi) $\Rightarrow$ (v) $\Rightarrow$ (iv)
 $\Rightarrow$ (iii) $\Rightarrow$ (ii) $\Rightarrow$ (i) are obvious.
The implication (i) $\Rightarrow$ (v) is a consequence of Proposition~\ref{local'}.

To show (v) $\Rightarrow$ (vi) set $m:=\dim M$ and consider any transitive family 
$\phi\colon (M,p)\times (\R^m,0) \to (M,p)$
of CR-automorphisms as in Definition~\ref{tr-family}.
Since $(M,p)$ is real-analytic, we may assume it is
embedded as a generic submanifold of $\C^N$.
Then the germ of a CR-map $\phi$ extends to a germ of a holomorphic map 
$\Phi\colon (\C^N,p)\times (\C^m,0) \to (\C^N,p)$.
Differentiating $\Phi$ in the second component in the direction of
the standard unit vectors in $\R^m$,
we obtain $m$ holomorphic vector fields $\C^N$
whose real parts are tangent to $M$.
Hence their restrictions to $M$ are infinitesimal CR-automorphisms.
Furthermore, by the assumption $\phi_*(T_0 \R^m)= T_pM$ in Definition~\ref{tr-family},
the values of these vector fields at $p$ span $T_p M$. 
This proves (vi).

To show (vi) $\Rightarrow$ (vii), we note that, by Proposition~\ref{local'},
$(M,p)$ is CR-equivalent to $(\2M\times\C^{N_2}\times\R^{N_3},0)$ with suitable 
$N_2$ and $N_3$ such that $(\2M,0)$ satisfies (v) and hence also (vi)
by the argument just before. Furthermore, $(\2M,0)$ is
both finitely nondegenerate and of finite type in view of Proposition~\ref{local'}~(ii).
Then a result by M.S.~Baouendi, P.~Ebenfelt and L.P.~Rothschild \cite{BERcomm}
implies that the Lie algebra of all germs of infinitesimal CR-automorphisms of $(\2M,0)$
is finite-dimensional. Since $(\2M,0)$ satisfies (vi), this Lie algebra must span $T_pM$.
Adding constant vector fields in the directions of $\C^{N_2}$ and $\R^{N_3}$
to this algebra, we easily conclude that also $(\2M\times\C^{N_2}\times\R^{N_3},0)$
satisfies (vii). Since the latter germ is CR-equivalent to $(M,p)$,
we also have (vii) for $(M,p)$.

Finally, given a finite-dimensional Lie algebra $\frak g$ as in (vii),
let $G$ be the corresponding connected and simply connected Lie group.
Then it is a well-known fact (Lie's Second Fundamental Theorem)
that $\frak g$ induces a local action of $G$ in a neighborhood of $p$ in $M$
such that the transformations by elements of $G$ correspond to local flows
of the vector fields from $\frak g$.
Since $\frak g$ consists of infinitesimal CR-automorphisms,
we conclude that the action obtained is by CR-automorphisms as desired.
The fact that the action of $G$ is transitive easily follows from 
the assumption in (vii) that $\frak g$ spans $T_pM$. This proves (viii),
completing the proof of Theorem~\ref{local}.

To prove Theorem~\ref{global}, we first note that
it follows directly from Theorem~\ref{local} that
conditions (iv) -- (vii) in Theorem~\ref{global} are equivalent.
Furthermore, by the equivalence of (v) and (iv) in Theorem~\ref{local},
it follows that (iv) in Theorem~\ref{global} implies
that every $p$ has a neighborhood $U(p)$ in $M$ such that
$(M,q)$ is CR-equivalent to $(M,p)$ for every $q\in U(p)$.
Since $M$ is connected, it is easy to see that its
germs at any two points are CR-equivalent, proving (iii).
Hence we have the implication (iv) $\Rightarrow$ (iii)
and the implications (iii) $\Rightarrow$ (ii) $\Rightarrow$ (i) are obvious.
Finally, applying again Theorem~\ref{local}, we see that (i) in Theorem~\ref{global}
implies (iv) there. Hence all conditions in Theorem~\ref{global}
are equivalent and the proof is complete.

\end{document}